# DEFAULT PRIORS FOR GAUSSIAN PROCESSES[1]

By Rui Paulo

*National Institute of Statistical Sciences and Statistical and Applied Mathematical Sciences Institute*

Motivated by the statistical evaluation of complex computer models, we deal with the issue of objective prior specification for the parameters of Gaussian processes. In particular, we derive the Jeffreys-rule, independence Jeffreys and reference priors for this situation, and prove that the resulting posterior distributions are proper under a quite general set of conditions. A proper flat prior strategy, based on maximum likelihood estimates, is also considered, and all priors are then compared on the grounds of the frequentist properties of the ensuing Bayesian procedures. Computational issues are also addressed in the paper, and we illustrate the proposed solutions by means of an example taken from the field of complex computer model validation.

**1. Introduction.** In this paper we address the problem of specifying objective priors for the parameters of general Gaussian processes. We derive formulas for the Jeffreys-rule prior, for an independence Jeffreys prior and for a reference prior on the parameters involved in the parametric specification of the mean and covariance functions of a Gaussian process. The mean is assumed to be a $q$-dimensional linear model on location-dependent covariates, and the correlation function involves an $r$-dimensional vector of unknown parameters. The resulting posteriors are shown to be proper under a more restrictive scenario. We also address computational issues, and in particular we devise a sampling scheme to draw from the resulting posteriors that requires very little input from the user. A method aimed at

Received December 2003; revised April 2004.

[1]Supported by NSF Grants DMS-00-73952 at the National Institute of Statistical Sciences and DMS-01-03265 at the Statistical and Applied Mathematical Sciences Institute. This research formed part of the author's Ph.D. dissertation at Duke University, where he was supported in part by a Ph.D. Fellowship awarded by Fundação para a Ciência e a Tecnologia, Portugal, with reference PRAXIS XXI/BD/15703/98.

producing proper "diffuse" priors, based on maximum likelihood estimates, is also described, and we present the results of a simulation study designed to compare the frequentist properties of all the Bayesian methods described in the paper. An example based on real data is used to illustrate some of the proposed solutions.

The motivation for considering the problems addressed in this paper is both theoretical and practical. From the applied point of view, these results are of interest in general for the field of spatial statistics, but especially for the analysis and validation of complex computer models. Indeed, one prominent approach to this problem involves fitting a Gaussian process to the computer model output, and a separable correlation function involving several parameters, typically a multidimensional power exponential, is frequently assumed—see, for example, Sacks, Welch, Mitchell and Wynn (1989), Kennedy and O'Hagan (2000, 2001) and Bayarri et al. (2002). The computer models are often computationally extremely demanding and this is a way of providing a cheaper surrogate that can be used in the design of computer experiments, optimization problems, prediction and uncertainty analysis, calibration and validation of the model. In the Bayesian approach, one must specify prior distributions for the parameters of the Gaussian processes. Typically, little or no prior information about the parameters is available, and their interpretation is not always straightforward, so that automatic or default procedures are sought. The need for default specification of priors, and also for the development of computational schemes for the resulting posteriors, is thus considerable in this area.

From a more theoretical perspective, this article is also relevant. Berger, De Oliveira and Sansó (2001) consider objective Bayesian analysis of Gaussian spatial processes with a quite general correlation structure, but the study is restricted to the situation where only the (one-dimensional) range parameter is considered to be unknown. The original motivation for their paper was the observation that commonly prescribed default priors could fail to yield proper posteriors. A more in-depth study of the problem revealed very interesting and unusual facts. In the presence of an unknown mean level for the Gaussian process, the integrated likelihood for the parameters governing the correlation structure is typically bounded away from zero, which explains the difficulty with posterior propriety. Also, the independence Jeffreys prior (assuming the parameters in the mean level are a priori independent of the ones involved in the covariance), which is often prescribed, fails to yield a proper posterior when the mean function includes an unknown constant level. The usual algorithm for the reference prior of Bernardo (1979) and Berger and Bernardo (1992), in which asymptotic marginalization is used, also fails to produce a proper posterior—this is the first known example in which exact marginalization is required to achieve posterior propriety. The authors end up recommending this "exact" reference prior.

The next logical step is to investigate whether these results still hold in higher dimensions, which is precisely one of the aims of the present article. The answer is no, in the sense that there are no surprises in terms of posterior propriety, and we describe in detail the reasons for that.

The paper is organized as follows. Section 2 sets up some notation and establishes formulas for the objective priors that are valid in a very general setting. The next section begins by describing the scenario where analytical results have been achieved, and in the sequel we analyze the behavior of the integrated likelihood and of the priors. Finally, conditions ensuring posterior propriety are determined.

Section 4 addresses computational aspects of the problem, and in particular we describe a Markov chain algorithm to sample from the posterior that requires very little input from the user. This algorithm involves computing maximum likelihood estimates and the Fisher information matrix. Taking advantage of the availability of these quantities, we propose also an empirical Bayes approach to the problem that aims at reproducing the practice of placing proper flat priors on the parameters. This section ends with an example, using real data, that illustrates some of the proposed solutions.

Section 5 presents the results of a simulation study designed to compare the frequentist properties of the Bayesian procedures proposed in the paper, and ends with some final recommendations.

In the Appendix we present the proofs of the various results that are described in the body of the paper.

**2. Notation and the objective priors.** Let us consider the following rather general situation: $Y(\cdot)$ is a Gaussian process on $\mathcal{S} \subset \mathbb{R}^p$ with mean and covariance functions given, respectively, by

$$\mathbb{E} Y(\mathbf{x}) \equiv \boldsymbol{\Psi}(\mathbf{x})' \boldsymbol{\theta} \quad \text{and} \quad \text{Cov}\,(Y(\mathbf{x}), Y(\mathbf{x}^\star)) \equiv \sigma^2 c(\mathbf{x}, \mathbf{x}^\star | \boldsymbol{\xi}),$$

where $c(\cdot, \cdot)$ is the correlation function, $\sigma^2$ is the variance and $\boldsymbol{\xi}$ is an $r$-dimensional vector of unknown positive parameters. The vector $\boldsymbol{\Psi}(\mathbf{x}) = (\psi_1(\mathbf{x}), \ldots, \psi_q(\mathbf{x}))'$ is a $q$-vector of location-dependent covariates. Define $\boldsymbol{\eta} = (\sigma^2, \boldsymbol{\theta}', \boldsymbol{\xi}')'$.

The stochastic process $Y(\cdot)$ is observed at locations $S = \{\mathbf{x}_1, \ldots, \mathbf{x}_n\}$ and hence the resulting random vector, $\mathbf{Y} = (Y(\mathbf{x}_1), \ldots, Y(\mathbf{x}_n))'$, satisfies

(2.1) $$\mathbf{Y} | \boldsymbol{\eta} \sim \mathsf{N}(\mathbf{X}\boldsymbol{\theta}, \sigma^2 \boldsymbol{\Sigma})$$

where $\boldsymbol{\Sigma} \equiv \boldsymbol{\Sigma}(\boldsymbol{\xi}) = [c(\mathbf{x}_i, \mathbf{x}_j | \boldsymbol{\xi})]_{ij}$ and

(2.2) $$\mathbf{X} = (\boldsymbol{\Psi}(\mathbf{x}_1)' \cdots \boldsymbol{\Psi}(\mathbf{x}_n)')'.$$

As a result, the associated likelihood function based on the observed data $\mathbf{y}$ is given by

(2.3) $$L(\boldsymbol{\eta} | \mathbf{y}) \propto (\sigma^2)^{-n/2} |\boldsymbol{\Sigma}|^{-1/2} \exp\left\{-\frac{1}{2\sigma^2}(\mathbf{y} - \mathbf{X}\boldsymbol{\theta})' \boldsymbol{\Sigma}^{-1} (\mathbf{y} - \mathbf{X}\boldsymbol{\theta})\right\}.$$

The priors on $\boldsymbol{\eta}$ that we will consider are of the form

$$\pi(\boldsymbol{\eta}) \propto \frac{\pi(\boldsymbol{\xi})}{(\sigma^2)^a} \tag{2.4}$$

for different $\pi(\boldsymbol{\xi})$ and $a$. Indeed, we will in the next two propositions show that this is the case for the Jeffreys-rule prior, for an independence Jeffreys prior and for a reference prior. This general prior follows a familiar form: flat on the parameters that specify the mean and usual forms for the variance $\sigma^2$.

Following Berger, De Oliveira and Sansó (2001), in order to derive the reference prior we specify the parameter of interest to be $(\sigma^2, \boldsymbol{\xi})$ and consider $\boldsymbol{\theta}$ to be the nuisance parameter. This corresponds in the reference prior algorithm to factoring the prior as

$$\pi^R(\boldsymbol{\eta}) = \pi^R(\boldsymbol{\theta}|\sigma^2, \boldsymbol{\xi})\pi^R(\sigma^2, \boldsymbol{\xi})$$

and selecting $\pi^R(\boldsymbol{\theta}|\sigma^2, \boldsymbol{\xi}) \propto 1$, because that is the Jeffreys-rule prior for the model at hand when $\sigma^2$ and $\boldsymbol{\xi}$ are considered known. Next, $\pi^R(\sigma^2, \boldsymbol{\xi})$ is calculated as the reference prior but for the marginal experiment defined by the integrated likelihood with respect to $\pi^R(\boldsymbol{\theta})$. It is this marginalization step that usually is carried out in an asymptotic fashion; Berger, De Oliveira and Sansó (2001) recommend that, for statistical models with a complicated covariance structure, the exact marginalization should become standard practice.

Before stating the formula for this reference prior, let us define $\dot{\boldsymbol{\Sigma}}^k$ as the matrix that results from $\boldsymbol{\Sigma}$ by differentiating each of its components with respect to $\xi_k$, the $k$th component of $\boldsymbol{\xi}$. Also, note that the integrated likelihood with respect to $\pi^R(\boldsymbol{\theta}) \propto 1$ is in this case

$$L^I(\sigma^2, \boldsymbol{\xi}|\mathbf{y}) = \int_{\mathbb{R}^q} L(\boldsymbol{\eta}|\mathbf{y})\pi^R(\boldsymbol{\theta})\,d\boldsymbol{\theta}$$
$$\propto (\sigma^2)^{-(n-q)/2}|\boldsymbol{\Sigma}|^{-1/2}|\mathbf{X}'\boldsymbol{\Sigma}^{-1}\mathbf{X}|^{-1/2}\exp\left\{-\frac{S^2_{\boldsymbol{\xi}}}{2\sigma^2}\right\}, \tag{2.5}$$

where $S^2_{\boldsymbol{\xi}} = \mathbf{y}'\mathbf{Q}\mathbf{y}$, $\mathbf{Q} = \boldsymbol{\Sigma}^{-1}\mathbf{P}$ and $\mathbf{P} = \mathbf{I} - \mathbf{X}(\mathbf{X}'\boldsymbol{\Sigma}^{-1}\mathbf{X})^{-1}\mathbf{X}'\boldsymbol{\Sigma}^{-1}$. Berger, De Oliveira and Sansó (2001), resorting to a result by Harville (1974), point out that there is a particular transformation of the data which has sampling distribution proportional to (2.5), and hence it is legitimate to compute the associated Jeffreys-rule prior.

PROPOSITION 2.1. *The reference prior $\pi^R(\boldsymbol{\eta})$ is of the form* (2.4) *with*

$$a = 1 \quad \text{and} \quad \pi^R(\boldsymbol{\xi}) \propto |I_R(\boldsymbol{\xi})|^{1/2}, \tag{2.6}$$

*where, with* $\mathbf{W}_k = \dot{\boldsymbol{\Sigma}}^k \mathbf{Q}$, $k = 1, \ldots, r$,

$$(2.7) \quad I_R(\boldsymbol{\xi}) = \begin{pmatrix} n-q & \operatorname{tr} \mathbf{W}_1 & \operatorname{tr} \mathbf{W}_2 & \cdots & \operatorname{tr} \mathbf{W}_r \\ & \operatorname{tr} \mathbf{W}_1^2 & \operatorname{tr} \mathbf{W}_1 \mathbf{W}_2 & \cdots & \operatorname{tr} \mathbf{W}_1 \mathbf{W}_r \\ & & \ddots & \cdots & \vdots \\ & & & & \operatorname{tr} \mathbf{W}_r^2 \end{pmatrix}.$$

For the proof see Appendix A.0.2.

In the next proposition, we present the formulas for the two Jeffreys-type priors we will consider.

PROPOSITION 2.2. *The independence Jeffreys prior,* $\pi^{J1}$, *obtained by assuming* $\boldsymbol{\theta}$ *and* $(\sigma^2, \boldsymbol{\xi})$ *a priori independent, and the Jeffreys-rule prior,* $\pi^{J2}$, *are of the form* (2.4) *with, respectively,*

$$(2.8) \quad a = 1 \quad \text{and} \quad \pi^{J1}(\boldsymbol{\xi}) \propto |I_J(\boldsymbol{\xi})|^{1/2}$$

*and*

$$(2.9) \quad a = 1 + \frac{q}{2} \quad \text{and} \quad \pi^{J2}(\boldsymbol{\xi}) \propto |\mathbf{X}'\boldsymbol{\Sigma}^{-1}\mathbf{X}|^{1/2} \pi^{J1}(\boldsymbol{\xi}),$$

*where* $\mathbf{U}_k = \dot{\boldsymbol{\Sigma}}^k \boldsymbol{\Sigma}^{-1}$ *and*

$$(2.10) \quad I_J(\boldsymbol{\xi}) = \begin{pmatrix} n & \operatorname{tr} \mathbf{U}_1 & \operatorname{tr} \mathbf{U}_2 & \cdots & \operatorname{tr} \mathbf{U}_r \\ & \operatorname{tr} \mathbf{U}_1^2 & \operatorname{tr} \mathbf{U}_1 \mathbf{U}_2 & \cdots & \operatorname{tr} \mathbf{U}_1 \mathbf{U}_r \\ & & \ddots & \cdots & \vdots \\ & & & & \operatorname{tr} \mathbf{U}_r^2 \end{pmatrix}.$$

For the proof see Appendix A.0.2.

For priors of the form (2.4), it is possible to integrate explicitly the product of the likelihood and the prior over $(\sigma^2, \boldsymbol{\theta})$. Indeed, standard calculations yield [as long as $a > 1 - (n-q)/2$]

$$\int_{\mathbb{R}^q \times \mathbb{R}_+} L(\boldsymbol{\eta}|\mathbf{y}) \pi(\boldsymbol{\eta}) \, d\boldsymbol{\theta} \, d\sigma^2 = L^I(\boldsymbol{\xi}|\mathbf{y}) \pi(\boldsymbol{\xi}),$$

where

$$(2.11) \quad L^I(\boldsymbol{\xi}|\mathbf{y}) \propto |\boldsymbol{\Sigma}|^{-1/2} |\mathbf{X}'\boldsymbol{\Sigma}^{-1}\mathbf{X}|^{-1/2} (S_{\boldsymbol{\xi}}^2)^{-((n-q)/2 + a - 1)}.$$

It is clear that the posterior associated with the prior (2.4) is proper if and only if $0 < \int_\Omega L^I(\boldsymbol{\xi}|\mathbf{y}) \pi(\boldsymbol{\xi}) \, d\boldsymbol{\xi} < \infty$, where $\Omega \subset \mathbb{R}^r$ is the parametric space of $\boldsymbol{\xi}$.

It appears to be difficult to study analytically the properties of both the integrated likelihood $L^I(\boldsymbol{\xi}|\mathbf{y})$ and the function $\pi(\boldsymbol{\xi})$—which we will often refer to, in an abuse of terminology, as the marginal prior for $\boldsymbol{\xi}$—in such a general setting. In the next section we will make the problem more amenable to analytical treatment by introducing a number of assumptions defining a more restrictive scenario.

**3. Posterior propriety.**  In this section we show that the objective priors yield proper posteriors for an important special case. The proof is constructed as follows. In Section 3.1 we introduce the notion of a separable correlation function and the restrictions we assume on the linear model specifying the mean function. Next, under added assumptions we describe the behavior of the integrated likelihood, whereas in Section 3.3 we achieve the same goal but for each of the objective priors determined in Section 2. Putting these results together, we are able to prove the theorem stated in Section 3.4.

3.1. *Separable correlation functions.*  Let $p \geq 2$ and consider the partition of a general element of $\mathbf{x} \in \mathcal{S} \subset \mathbb{R}^p$ given by $\mathbf{x} = (\mathbf{x}_1, \mathbf{x}_2, \ldots, \mathbf{x}_r)$, where $r \leq p$ and each subvector $\mathbf{x}_k$, $k = 1, \ldots, r$, is of dimension $p_k$ with $\sum_{k=1}^{r} p_k = p$.

Denote by $d(\mathbf{x}, \mathbf{y})$ the Euclidean distance between two vectors. It is a simple consequence of Bochner's theorem [cf. Cressie (1993)] that if $c_k(\mathbf{x}_k, \mathbf{x}_k^\star) \equiv c_k(d(\mathbf{x}_k, \mathbf{x}_k^\star))$, $k = 1, \ldots, r$, are isotropic correlation functions in $\mathbb{R}^{p_k}$, then $c(\mathbf{x}, \mathbf{x}^\star) = \prod_{k=1}^{r} c_k(d(\mathbf{x}_k, \mathbf{x}_k^\star))$ is a valid correlation function in $\mathbb{R}^p$. Such correlation functions are called *separable*, or to be more precise, partially separable (the fully separable case corresponds to the choice $r = p$, $p_k \equiv 1$).

If a separable correlation function is used and if furthermore the set of locations at which the process is observed forms a suitable Cartesian product, it is easy to see that the correlation matrix of the data is the Kronecker product of the individual correlation matrices associated with each dimension $\mathbf{x}_k$. (Recall that the Kronecker product of two matrices, $\mathbf{A} = [a_{ij}]$ and $\mathbf{B}$, is defined by $\mathbf{A} \otimes \mathbf{B} = [a_{ij} \mathbf{B}]$.) It is essentially in this setting that we will study the analytical properties of the integrated likelihood and priors, along with establishing sufficient conditions for posterior propriety. We next make these statements precise.

We will henceforth assume that the following conditions hold.

ASSUMPTION $\mathcal{A}1$.  Separability of the correlation function: with $\mathbf{x} = (\mathbf{x}_1, \mathbf{x}_2, \ldots, \mathbf{x}_r)$, where $p \geq r \geq 2$,

$$c(\mathbf{x}, \mathbf{x}^\star | \boldsymbol{\xi}) = \prod_{k=1}^{r} \rho(d(\mathbf{x}_k - \mathbf{x}_k^\star); \xi_k),$$

where $\rho$ is a valid (isotropic) correlation function.

ASSUMPTION $\mathcal{A}2$.  Cartesian product of the design set:

$$S = S_1 \times S_2 \times \cdots \times S_r,$$

where $S_k = \{\mathbf{x}_{1k}, \ldots, \mathbf{x}_{n_k, k}\} \subset \mathbb{R}^{p_k}$, $\#S_k = n_k$, so that $\#S \equiv n = \prod_{k=1}^{r} n_k$.

ASSUMPTION $\mathcal{A}3$. *Mean structure*: $q = 1$ and

$$\mathbf{X} = \mathbf{X}_1 \otimes \mathbf{X}_2 \otimes \cdots \otimes \mathbf{X}_r$$

with each $\mathbf{X}_k$ of dimension $n_k \times 1$.

A few remarks are in order. First, one should think about the parameter $\xi_k$ as representing the range parameter of the associated correlation structure; if other parameters are present, for example, roughness parameters, those will be assumed known.

In terms of the mean structure, within the framework defined by Assumption $\mathcal{A}3$ we can find, for instance, the unknown constant mean case, that is, $\mathbf{\Psi}(\mathbf{x}) \equiv 1$ so that $\mathbf{X} \equiv \mathbf{1}_n = \mathbf{1}_{n_1} \otimes \cdots \otimes \mathbf{1}_{n_p}$. This is in turn an instance of the more general situation described by $\mathbf{\Psi}(\mathbf{x}) \equiv \prod_{k \in A} \psi_k(\mathbf{x}_k)$ where $A \subset \{1, \ldots, n\}$. In this circumstance, $\mathbf{X}_k = (\psi_k(\mathbf{x}_{ik}),\ i = 1, \ldots, n_k)'$ if $k \in A$; otherwise $\mathbf{X}_k = \mathbf{1}_{n_k}$.

As we alluded to before, Assumptions $\mathcal{A}1$ and $\mathcal{A}2$ jointly allow for a very convenient Kronecker product expression for the correlation matrix of the data. To make that clear and to set up some notation as well, we have

$$(3.1) \qquad \mathbf{\Sigma} = \mathbf{\Sigma}_1 \otimes \mathbf{\Sigma}_2 \otimes \cdots \otimes \mathbf{\Sigma}_r \equiv \bigotimes_{k=1}^{r} \Sigma_k,$$

where $\mathbf{\Sigma}_k = [\rho(d(\mathbf{x}_{ik} - \mathbf{x}_{jk}); \xi_k)]_{i,j=1,\ldots,n_k},\ k = 1, \ldots, r$, are the correlation matrices associated with each of the separated dimensions $\mathbf{x}_k$. Note that each of these matrices is of dimension $n_k \times n_k$. This fact is further explored in Section 4. For simplicity we will from here on use the shorthand exemplified above: we represent by $\bigotimes_{k=1}^{r} \mathbf{A}_k$ the matrix $\mathbf{A}_1 \otimes \cdots \otimes \mathbf{A}_r$.

The $\rho$ functions we will consider are quite general, but we will focus particular attention on the families of correlation functions that we list in Table 1 above for future reference. For more details and additional references, consult Cressie (1993).

3.2. *Behavior of the integrated likelihood.* In this section we will study properties of the integrated likelihood. Nevertheless, the next two lemmas are key also to the series of results that will follow concerning the analytical behavior of the priors.

LEMMA 3.1. *Under Assumptions $\mathcal{A}1$–$\mathcal{A}3$, we have*

$$(3.2) \qquad \mathbf{Q} = \mathbf{\Sigma}^{-1} \mathbf{P} = \bigotimes_{k=1}^{r} \mathbf{\Sigma}_k^{-1} - \bigotimes_{k=1}^{r} \mathbf{\Phi}_k,$$

*where* $\mathbf{\Phi}_k = \mathbf{\Sigma}_k^{-1} \mathbf{X}_k (\mathbf{X}_k' \mathbf{\Sigma}_k^{-1} \mathbf{X}_k)^{-1} \mathbf{X}_k' \mathbf{\Sigma}_k^{-1}$, *and*

$$(3.3) \qquad |\mathbf{X}' \mathbf{\Sigma}^{-1} \mathbf{X}| = \prod_{k=1}^{r} \mathbf{X}_k' \mathbf{\Sigma}_k^{-1} \mathbf{X}_k.$$

TABLE 1
*Families of correlation functions*

Spherical—
$$\rho(d;\xi) = [1 - \tfrac{3}{2}d\xi + \tfrac{1}{2}(d\xi)^3]I\{d\xi \leq 1\};\ \xi > 0.$$
Power exponential—
$$\rho(d;\xi) = \exp[-(d\xi)^\alpha];\ \xi > 0,\ \alpha \in (0, 2].$$
Rational quadratic—
$$\rho(d;\xi) = [1 + (d\xi)^2]^{-\alpha};\ \xi > 0,\ \alpha > 0.$$
Matérn— $\mathcal{K}_\alpha(\cdot)$ is the modified Bessel function of the second kind and order $\alpha$:
$$\rho(d;\xi) = \frac{1}{2^{\alpha-1}\Gamma(\alpha)}(d\xi)^\alpha \mathcal{K}_\alpha(d\xi);\ \xi > 0,\ \alpha > 0.$$

$\xi$ denotes a range parameter; $\alpha$ refers to a roughness parameter; $d$ represents the Euclidean distance between two locations.

For the proof see Appendix A.0.4.

LEMMA 3.2. *Under Assumptions $\mathcal{A}1$–$\mathcal{A}3$, suppose $\rho(d;\xi)$ is a continuous function of $\xi > 0$ for every $d \geq 0$ such that:*

(i) *$\rho(d;\xi) = \rho^0(d\xi)$, where $\rho^0(\cdot)$ is a correlation function satisfying $\lim_{u\to\infty}\rho^0(u) = 0$;*
(ii) *as $\xi_k \to 0$, each of the correlation matrices $\mathbf{\Sigma}_k$ satisfies*

$$(3.4) \qquad \mathbf{\Sigma}_k = \mathbf{1}_{n_k}\mathbf{1}'_{n_k} + \nu(\xi_k)(\mathbf{D}_k + o(1)), \qquad k = 1,\ldots,r,$$

*for some continuous nonnegative function $\nu(\cdot)$ and fixed* nonsingular *matrix $\mathbf{D}_k$;*

(iii) *above, and for $k = 1, \ldots, r$, $\mathbf{D}_k$ should satisfy $\mathbf{1}'_{n_k}\mathbf{D}_k^{-1}\mathbf{1}_{n_k} \neq 0$, $\mathbf{X}'_k\mathbf{D}_k^{-1}\mathbf{X}_k \neq 0$ and, if $\mathbf{X}_k \neq \mathbf{1}$,*

$$(3.5) \qquad \mathbf{1}'_{n_k}\mathbf{D}_k^{-1}\mathbf{1}_{n_k} \neq \mathbf{1}'_{n_k}\mathbf{D}_k^{-1}\mathbf{X}_k(\mathbf{X}'_k\mathbf{D}_k^{-1}\mathbf{X}_k)^{-1}\mathbf{X}'_k\mathbf{D}_k^{-1}\mathbf{1}_{n_k}.$$

*Define the quantities*

$$(3.6) \qquad \mathbf{F}_k = \frac{\mathbf{D}_k^{-1}\mathbf{1}_{n_k}\mathbf{1}'_{n_k}\mathbf{D}_k^{-1}}{(\mathbf{1}'_{n_k}\mathbf{D}_k^{-1}\mathbf{1}_{n_k})^2},$$

$$(3.7) \qquad \mathbf{G}_k = \mathbf{D}_k^{-1} - \frac{\mathbf{D}_k^{-1}\mathbf{1}_{n_k}\mathbf{1}'_{n_k}\mathbf{D}_k^{-1}}{\mathbf{1}'_{n_k}\mathbf{D}_k^{-1}\mathbf{1}_{n_k}},$$

$$(3.8) \qquad \mathbf{H}_k = \frac{\mathbf{G}_k\mathbf{X}_k\mathbf{X}'_k\mathbf{G}_k}{\mathbf{X}'_k\mathbf{G}_k\mathbf{X}_k}.$$

*In this setting we have that, as $\xi_k \to 0$, $k = 1, \ldots, r$,*

$$(3.9) \qquad \mathbf{\Sigma}_k^{-1} = \mathbf{G}_k(1 + o(1))/\nu(\xi_k),$$

$$(3.10) \qquad |\mathbf{\Sigma}_k| = [\nu(\xi_k)]^{n_k - 1}|\mathbf{D}_k|(\mathbf{1}'_{n_k}\mathbf{D}_k\mathbf{1}_{n_k})(1 + o(1)),$$

$$(3.11) \quad \mathbf{X}_k' \boldsymbol{\Sigma}_k^{-1} \mathbf{X}_k = \begin{cases} 1 + o(1), & \text{if } \mathbf{X}_k = \mathbf{1}, \\ \left[ \mathbf{X}_k' \mathbf{D}_k^{-1} \mathbf{X}_k - \dfrac{\mathbf{1}' \mathbf{D}_k^{-1} \mathbf{X}_k \mathbf{X}_k' \mathbf{D}_k^{-1} \mathbf{1}}{\mathbf{1}' \mathbf{D}_k^{-1} \mathbf{1}} \right] (1 + o(1))/\nu(\xi_k), \\ \qquad \text{otherwise}, \end{cases}$$

$$(3.12) \quad \boldsymbol{\Phi}_k = \begin{cases} \mathbf{F}_k(1 + o(1)), & \text{if } \mathbf{X}_k = \mathbf{1}, \\ \mathbf{H}_k(1 + o(1))/\nu(\xi_k), & \text{otherwise.} \end{cases}$$

For the proof see Appendix A.0.4.

The following result describes the behavior of the integrated likelihood (2.11) in the present setting and under the assumptions of the previous lemma.

PROPOSITION 3.3. *Under the conditions of Lemma 3.2, and assuming a prior of the form* (2.4), $L^I(\boldsymbol{\xi}|\mathbf{y})$ *is a continuous function of $\boldsymbol{\xi}$ given by the expression*

$$(3.13) \quad L^I(\boldsymbol{\xi}|\mathbf{y}) \propto \prod_{k=1}^r \{|\boldsymbol{\Sigma}_k|^{-n_{(k)}/2} (\mathbf{X}_k' \boldsymbol{\Sigma}_k^{-1} \mathbf{X}_k)^{-1/2}\} \times [\mathbf{y}' \mathbf{Q} \mathbf{y}]^{-(n-3+2a)/2},$$

*where $\mathbf{Q}$ is given by* (3.2) *and* $n_{(k)} \equiv \prod_{i \neq k} n_i$.

*Additionally, let* $A = \{k \in \{1, \ldots, r\} : \mathbf{X}_k \neq \mathbf{1}\}$ *and* $\varnothing \neq B \subset \{1, \ldots, r\}$. *Denote by $\bar{E}$ the complement of a set $E$. Then,*

(a) *as $\xi_k \to 0$, $k \in B$ and $\xi_k \to \infty$, $k \in \bar{B}$,*

$$(3.14) \quad L^I(\boldsymbol{\xi}|\mathbf{y}) \propto \prod_{k \in B} \nu(\xi_k)^{(n_k - 3 + 2a)/2} \prod_{k \in A \cap B} \nu(\xi_k)^{1/2} (1 + g(\boldsymbol{\xi}))$$

*with $g(\boldsymbol{\xi}) \to 0$;*

(b) *as $\xi_k \to 0$, $k \in B$ and $0 < \delta_k < \xi_k < M_k < \infty$, $k \in \bar{B}$,*

$$L^I(\boldsymbol{\xi}|\mathbf{y}) \propto \prod_{k \in B} \nu(\xi_k)^{(n_k - 3 + 2a)/2} \prod_{k \in A \cap B} \nu(\xi_k)^{1/2}$$
$$(3.15) \qquad \qquad \times (1 + g(\{\xi_k, k \in B\})) \times h(\{\xi_k, k \in \bar{B}\}),$$

*with $g(\cdot) \to 0$ and $h(\cdot)$ continuous on $\{\delta_k < \xi_k < M_k, k \in \bar{B}\}$.*

For the proof see Appendix A.0.4.

3.3. *Behavior of the priors.* In the next two results, properties of both Jeffreys priors and of the reference prior are described. They are based on a series of assumptions that we list below and form (except for Assumption $\mathcal{A}8$) a subset of those introduced in Berger, De Oliveira and Sansó (2001) in their study of the asymptotic properties of the priors. According to this paper, the families of correlation functions listed in Table 1 (virtually always) satisfy these properties [e.g., beware of the fact that for the power exponential with $\alpha = 2$ and more than three equally spaced points, it is not the case that $\mathbf{D}$ is invertible. It is if $\alpha \in (0,2)$]. Note also how these assumptions imply those of Lemma 3.2. Below we denote by $\dot{\boldsymbol{\Sigma}}_k$ the matrix of derivatives of $\boldsymbol{\Sigma}_k$ with respect to $\xi_k$, $k = 1, \ldots, r$.

*Assumptions.* Suppose that $\rho(d; \xi) = \rho^0(d\,\xi)$, where $\rho^0(\cdot)$ is a correlation function satisfying $\lim_{u \to \infty} \rho^0(u) = 0$, is a continuous function of $\xi > 0$ for any $d > 0$ and that, for $k = 1, \ldots, r$,

ASSUMPTION $\mathcal{A}4$. As $\xi_k \to 0$, there are fixed matrices $\mathbf{D}_k$, *nonsingular*, and $\mathbf{D}_k^\star$; differentiable functions $\nu(\cdot)$ ($> 0$) and $\omega(\cdot)$; and a matrix function $\mathbf{R}_k(\cdot)$ that is differentiable as well, such that

$$(3.16) \qquad \boldsymbol{\Sigma}_k = \mathbf{1}_{n_k} \mathbf{1}'_{n_k} + \nu(\xi_k) \mathbf{D}_k + \omega(\xi_k) \mathbf{D}_k^\star + \mathbf{R}_k(\xi_k),$$

$$(3.17) \qquad \dot{\boldsymbol{\Sigma}}_k = \nu'(\xi_k) \mathbf{D}_k + \omega'(\xi_k) \mathbf{D}_k^\star + \frac{\partial}{\partial \xi_k} \mathbf{R}_k(\xi_k).$$

ASSUMPTION $\mathcal{A}5$. The functions $\nu$, $\omega$ and $\mathbf{R}_k$ further satisfy, as $\xi_k \to 0$ and with $\|[a_{ij}]\|_\infty = \max_{ij}\{|a_{ij}|\}$,

$$\frac{\omega(\xi_k)}{\nu(\xi_k)} \to 0, \qquad \frac{\omega'(\xi_k)}{\nu'(\xi_k)} \to 0,$$

$$\frac{\|\mathbf{R}_k(\xi_k)\|_\infty}{\nu(\xi_k)} \to 0, \qquad \frac{\|(\partial/(\partial \xi_k))\mathbf{R}_k(\xi_k)\|_\infty}{\nu'(\xi_k)} \to 0.$$

ASSUMPTION $\mathcal{A}6$. Above $\mathbf{D}_k$ satisfies $\mathbf{1}'_{n_k} \mathbf{D}_k^{-1} \mathbf{1}_{n_k} \neq 0$, $\mathbf{X}'_k \mathbf{D}_k^{-1} \mathbf{X}_k \neq 0$ and, if $\mathbf{X}_k \neq \mathbf{1}$,

$$\mathbf{1}'_{n_k} \mathbf{D}_k^{-1} \mathbf{1}_{n_k} \neq \mathbf{1}'_{n_k} \mathbf{D}_k^{-1} \mathbf{X}_k (\mathbf{X}'_k \mathbf{D}_k^{-1} \mathbf{X}_k)^{-1} \mathbf{X}'_k \mathbf{D}_k^{-1} \mathbf{1}_{n_k}.$$

ASSUMPTION $\mathcal{A}7$. $[\mathrm{tr}(\dot{\boldsymbol{\Sigma}}_k)^2]^{1/2}$, as a function of $\xi_k$, is integrable at infinity.

ASSUMPTION $\mathcal{A}8$. $n_k \geq 2$.

We start with the reference prior and then proceed to the Jeffreys-type priors.

PROPOSITION 3.4. *For the reference prior* (2.6), *and under Assumptions* $\mathcal{A}1$–$\mathcal{A}3$, *there are functions* $\pi_k^R(\xi_k)$, $k = 1, \ldots, r$, *such that*

$$(3.18) \qquad \pi^R(\boldsymbol{\xi}) \leq \prod_{k=1}^r \pi_k^R(\xi_k),$$

where, as $\xi_k \to 0$, $k = 1, \ldots, r$, and under the added Assumptions A4–A8,

$$\pi_k^R(\xi_k) \propto \left|\frac{\nu'(\xi_k)}{\nu(\xi_k)}\right|(1 + o(1)), \tag{3.19}$$

and $\pi_k^R$ is integrable at infinity.

For the proof see Appendix A.0.5.

PROPOSITION 3.5. *Under Assumptions A1–A3, for the independence Jeffreys prior $\pi^{J1}$ given by* (2.8), *and the Jeffreys-rule prior $\pi^{J2}$ given by* (2.9), *there are functions $\pi_k^{Ji}(\xi_k)$, $i = 1, 2$, $k = 1, \ldots, r$, such that*

$$\pi^{Ji}(\boldsymbol{\xi}) \leq \prod_{k=1}^r \pi_k^{Ji}(\xi_k), \tag{3.20}$$

*where, under added Assumptions A4–A8 and for $k = 1, \ldots, r$, we have as $\xi_k \to 0$,*

$$\pi_k^{Ji}(\boldsymbol{\xi}_k) \propto \frac{|\nu'(\xi_k)|}{|\nu(\xi_k)|^{\alpha_k}}(1 + o(1)) \tag{3.21}$$

*with*

$$\alpha_k = \begin{cases} 1, & \text{if } i = 1 \text{ or } (i = 2 \text{ and } \mathbf{X}_k = \mathbf{1}), \\ 3/2, & \text{if } i = 2 \text{ and } \mathbf{X}_k \neq \mathbf{1}. \end{cases} \tag{3.22}$$

*Also, $\pi_k^{Ji}$ is integrable at infinity.*

For the proof see Appendix A.0.5.

3.4. *Results on posterior propriety.* In this section we investigate the propriety of the formal posterior associated with a general prior of the form (2.4), in the scenario described throughout the previous sections. We end with a result that specifically addresses the case of both Jeffreys priors and the reference prior.

THEOREM 3.6. *Under the set of Assumptions A1–A8, the posterior associated with the general prior* (2.4)—*where $\pi(\boldsymbol{\xi})$ is either $\pi^R$, $\pi^{J1}$ or $\pi^{J2}$— is proper as long as*

$$a > 1/2. \tag{3.23}$$

For the proof see Appendix A.0.6.

If we restrict attention to the instances of the general prior that are of particular interest, then we have the following:

COROLLARY 3.7. *Under Assumptions $\mathcal{A}1$–$\mathcal{A}8$, the reference, the independence Jeffreys and the Jeffreys-rule priors yield proper posteriors.*

PROOF. Just recall that $a = 1$ in the case of the reference and independence Jeffreys prior, and that for the Jeffreys-rule $a = 1 + q/2$, where by assumption $q = 1$. □

3.5. *Discussion and generalizations.* It is interesting to understand what makes the multidimensional problem different from the unidimensional one, the case covered in Berger, De Oliveira and Sansó (2001). As we mentioned in the Introduction, the reason why posterior propriety is difficult to achieve in the unidimensional case has to do with the behavior of the integrated likelihood near the origin. For example, if $p = 1$ and $\mathbf{X} = \mathbf{1}$, we have, as $\xi \to 0$, $L^I(\xi|\mathbf{y}) = O([\nu(\xi)]^{1-a})$, which is independent of the sample size. Roughly speaking, the first and last factors of (2.11) produce powers of $\nu(\xi)$ depending on the sample size that essentially cancel.

The key formula in the multidimensional case is

$$\left| \bigotimes_{k=1}^{r} \boldsymbol{\Sigma}_k \right| = \prod_{k=1}^{r} |\boldsymbol{\Sigma}_k|^{n_{(k)}}$$

where $n_{(k)} = \prod_{i \neq k} n_i$. Considering (3.10), it is clear that, as $\xi_k \to 0$, $|\boldsymbol{\Sigma}|^{-1/2} \propto [\prod_{i \neq k}^{r} |\boldsymbol{\Sigma}_i|^{-n_{(i)}/2}]\nu(\xi_k)^{-n/2}\nu(\xi_k)^{n_{(k)}/2}(1 + o(1))$. The first factor involving $\nu(\xi_k)$ will again essentially cancel, but we still keep the second, which explains the different behaviors. For small $\xi_k$ the prior is not that relevant in determining posterior propriety in the multidimensional case, whereas in the unidimensional case it is of capital importance.

Another point of interest is the range of applicability of the results of this section. We would argue that Assumptions $\mathcal{A}2$ and $\mathcal{A}3$ are mathematical conveniences and, as we have stated before, Assumptions $\mathcal{A}4$–$\mathcal{A}7$ are virtually always satisfied in the context of the usual families of correlation functions. The (partial) separability Assumption $\mathcal{A}1$ is indeed restrictive, and more research is needed before one feels confident about using the priors derived in Section 2 in circumstances where Assumption $\mathcal{A}1$ is not valid.

Nevertheless, Assumption $\mathcal{A}1$ is not uncommon in practice. Apart from the field of computer model validation, where it is frequently assumed, various other examples can be found in geology applications, where the two lateral dimensions are separated from the vertical dimension, and in space-time models, where time is separated from the position vector [consider, e.g., Fuentes (2003), Short and Carlin (2003) and references therein].

We would also like to stress the theoretical interest of the results of this section. As we pointed out, Berger, De Oliveira and Sansó (2001) revealed very unusual and surprising facts unknown to the spatial statistics community at the time. There was no apparent reason to expect a different scenario in any sensible extension of the setting treated in that paper. The present article seems to indicate that, possibly, Berger, De Oliveira and Sansó (2001) deal with an exceptional case, and that, in general, the "standard" type of behavior is to be expected.

There are obvious potential extensions to the present work. One would like to treat the smoothness parameters as unknowns and also to consider variations on the ordering of the parameter vector in the reference prior. These are technically challenging tasks that constitute work in progress.

**4. Computation.** There is considerable need in the area of computer model validation for default prior specifications and also for efficient computational schemes. In this section we address several computing issues, and in particular that of Bayesian learning in the presence of the objective priors derived in this paper. We present an example, taken from the field of computer model validation, to illustrate some of the proposed solutions.

We now return to the general setting of Section 2, and in particular, unless otherwise noted, $\boldsymbol{\xi}$ represents a general $r$-dimensional vector involved in the parametric formulation of the correlation function.

4.1. *General remarks.* Evaluating any of the objective priors at one particular value of $\boldsymbol{\xi}$ is a computationally intensive task. The reference prior is the most computationally intensive of all three, followed by the independence Jeffreys prior and by the Jeffreys-rule prior. In order to convince oneself of that, it suffices to consider the calculations involved in computing each of the entries of the matrices (2.7) and (2.10). Since $\mathbf{W}_k = \mathbf{U}_k\,\mathbf{P}$, to compute the reference prior at a particular value of $\boldsymbol{\xi}$ one has to calculate the projection matrix $\mathbf{P}$, which is not necessary in any of the other priors. Additionally, each $\mathbf{W}_k$ requires one more matrix product (of two $n \times n$ matrices) than $\mathbf{U}_k$. The difference between the Jeffreys-rule and independence Jeffreys is only in the computation of $|\mathbf{X}'\boldsymbol{\Sigma}^{-1}\mathbf{X}|$. For similar reasons, it is clear that the likelihood function (2.3) is computationally less expensive than either of the integrated likelihoods given by (2.5) and (2.11). These observations are particularly relevant in the context of Markov chain Monte Carlo (MCMC) computations, addressed in Section 4.4, where the prior and the likelihood have to be evaluated a very large number of times.

One should note that some of the assumptions of Section 3 have considerable potential impact on the computational side of the problem. When satisfied, the separability Assumption $\mathcal{A}1$ paired with the Cartesian product structure of the design set Assumption $\mathcal{A}2$, which implies (3.1), can be exploited in order to tremendously simplify and speed up computation. Indeed, the most expensive (and possibly numerically unstable) calculations of this problem are certainly computing the inverse of the correlation matrix and its determinant. Facts 2 and 3 of Appendix A.0.1 essentially state that, in this case, we only have to compute the inverse and the determinant of each of the $p$ matrices $\boldsymbol{\Sigma}_k$ of dimension $n_k \times n_k$, and not for the $n \times n$ correlation matrix $\boldsymbol{\Sigma}$. Even the Cholesky decomposition of $\boldsymbol{\Sigma}$ can be obtained from that of the $\boldsymbol{\Sigma}_k$, as it is easy to see. This allows for much more freedom on the number of points at which one observes the stochastic process, as explored in Bayarri et al. (2002) when dealing with functional output of a computer model.

4.2. *Maximum likelihood estimates.* The calculation of maximum likelihood (ML) estimates turns out to be useful in several ways.

Indeed, in Section 4.3 we are going to use maximum likelihood estimates in devising a mechanism that produces proper flat priors. Also, these estimates can be useful in devising sampling schemes that require little input from the user, as we will see in Section 4.4. Last, it is often the case that there is not enough information in the data to learn about all the parameters involved in the correlation structure, especially about the parameters that control geometric properties of the process, the so-called roughness parameters. We have observed that as long as one fixes those parameters at "sensible" values, inference is fairly insensitive to the particular value chosen. Data-dependent choices, and in particular maximum likelihood estimates, have produced good practical results in our applied work.

The question of which kind of estimate one should compute is relevant, as here we have at our disposal explicit formulas for three distinct (integrated) likelihood functions.

There is considerable empirical evidence supporting the fact that, in general, maximum likelihood estimates derived from integrated likelihoods tend to be more stable and also more meaningful, so that one would in principle discard the idea of maximizing the joint likelihood (2.3). Also, these estimates are not as useful in the context of devising sampling schemes.

Maximizing the integrated likelihood $L^I(\sigma^2, \boldsymbol{\xi}|\mathbf{y})$ in (2.5) has some advantages over maximizing $L^I(\boldsymbol{\xi}|\mathbf{y})$ in (2.11). First of all, it gives rise to considerably simpler formulas for the gradient and the associated expected information matrix. Indeed, it is easy to see from inspection of the proof of Proposition 2.1 that

$$\frac{\partial}{\partial \sigma^2} \ln L^I(\sigma^2, \boldsymbol{\xi}) = \frac{1}{2}\left(\frac{1}{\sigma^4}\mathbf{y}'\mathbf{Q}\mathbf{y} - \frac{1}{\sigma^2}(n-q)\right),$$

$$\frac{\partial}{\partial \xi_k} \ln L^I(\sigma^2, \boldsymbol{\xi}) = \frac{1}{2}\left(\frac{1}{\sigma^2}\mathbf{y}'\mathbf{Q}\mathbf{W}_k\mathbf{y} - \operatorname{tr}\mathbf{W}_k\right), \qquad k=1,\ldots,r,$$

and that the associated Fisher information matrix is

$$(4.1) \quad I(\sigma^2, \boldsymbol{\xi}) = \frac{1}{2} \begin{pmatrix} \frac{n-q}{\sigma^4} & \frac{1}{\sigma^2}\operatorname{tr}\mathbf{W}_1 & \frac{1}{\sigma^2}\operatorname{tr}\mathbf{W}_2 & \cdots & \frac{1}{\sigma^2}\operatorname{tr}\mathbf{W}_r \\ & \operatorname{tr}\mathbf{W}_1^2 & \operatorname{tr}\mathbf{W}_1\mathbf{W}_2 & \cdots & \operatorname{tr}\mathbf{W}_1\mathbf{W}_r \\ & & \ddots & \cdots & \vdots \\ & & & & \operatorname{tr}\mathbf{W}_r^2 \end{pmatrix}.$$

This allows for a very simple optimization algorithm to be used in computing the associated estimate, namely Fisher's scoring method, a variant of the Newton–Raphson method that results from approximating the Hessian of the logarithm of integrated likelihood by its expected value. To be more precise, the $(s+1)$st iterate of the numerical method is given by

$$(\sigma^2)^{(s+1)} = \frac{\mathbf{y}'\mathbf{Q}\mathbf{y}}{n-q},$$

$$\boldsymbol{\xi}^{(s+1)} = \boldsymbol{\xi}^{(s)} + \lambda[I(\boldsymbol{\xi}^{(s)})]^{-1} \frac{\partial \ln L^I(\boldsymbol{\xi}, \sigma^2)}{\partial \boldsymbol{\xi}}\bigg|_{\boldsymbol{\xi}=\boldsymbol{\xi}^{(s)}, \sigma^2=(\sigma^2)^{(s)}},$$

where $I(\boldsymbol{\xi})$ results from (4.1) by dropping the first row and column. The quantity $\lambda$ is the step size of the algorithm.

There are some drawbacks to this simple numerical method. First, it is possible that, initially, some iterates happen to lie outside the parameter space. We avoid this problem by simply saying that $\xi_k^{(s+1)} = \xi_k^{(s)}$ whenever the $k$th component of $\boldsymbol{\xi}^{(s+1)}$ happens to not belong to the corresponding parameter space. Also, it is well known that Newton–Raphson-type methods are quite sensitive to the starting points, sometimes becoming trapped in local maxima and sometimes simply not converging. This unfortunately is not easy to solve. Some experimentation and tuning of $\lambda$ are required in order to assure convergence.

Getting an estimate of $\boldsymbol{\xi}$ by maximizing (2.11) is not as attractive because of the fact that the formulas are computationally more involved. For example, there does not seem to be a closed-form expression for the associated expected information matrix. Also, in the examples we have dealt with involving the power exponential family of correlation functions, we did not notice any improvements over the simpler method. The availability of an estimate of the variance will also be useful in the next section.

4.3. *Proper flat priors.* When trying to produce so-called noninformative Bayesian analyses, people often consider placing proper flat priors on the parameters as an alternative to using priors derived using formal methods. The basic idea is to specify priors that are relatively flat in the region of the parametric space where most of the posterior mass accumulates.

One way of reproducing this practice is by describing it as an empirical Bayes method. We do so by placing independent exponential priors centered at a multiple of the ML estimates on the precision (reciprocal of variance) and on the components of vector $\boldsymbol{\xi}$. The constant that multiplies the ML estimate in order to get the mean of the prior can be determined by experimentation, making sure that the effect of the prior on the posterior is relatively small, in other words, making sure that the prior is relatively flat in the region of the parametric space where the likelihood has most of its mass. For an example, consider Section 4.5.

In Section 5 we will compare the frequentist properties of this approach to objective Bayesian analysis to the use of the formal priors we derive in this paper. Apart from philosophical reasons, this comparison is also relevant because in the present case the empirical Bayes method is computationally much simpler than any of the alternatives.

4.4. *Sampling from the posterior.* No matter which prior specification one uses—any of the objective priors or the empirical Bayes method—one can sample exactly from the full conditional of $\boldsymbol{\theta}$ and of $\sigma^2$: these are, respectively, Gaussian and inverse-gamma.

Then, one has the choice between sampling from the marginal posterior $[\boldsymbol{\xi}|\mathbf{y}]$ or from the full conditional $[\boldsymbol{\xi}|\mathbf{y}, \sigma^2, \boldsymbol{\theta}]$. We adopt this latter strategy since the integrated likelihood is computationally more expensive than the full likelihood and, based on the experiments we conducted, there are no apparent gains in terms of mixing in adopting the other strategy.

At this point we have also the option of drawing $\boldsymbol{\xi}$ as a block or not. Notice that as long as one updates one of the components of $\boldsymbol{\xi}$ one has basically to recompute the likelihood, as it does not factor in any useful way. From that perspective, it is more efficient therefore to sample $\boldsymbol{\xi}$ as a block. Additionally, it is reasonable to expect some serial correlation between the $\xi_k$ in the MCMC output, and sampling $\boldsymbol{\xi}$ as a block would at least reduce it.

The question is, can we get a proposal, in a more or less automatic fashion, that allows for this? The ML estimate and the expected information matrix we described in Section 4.2 are useful in accomplishing this goal. Although the resulting sampling scheme is certainly not new and has been widely used, it seems to work well in the examples we have examined so far.

The idea is to consider a reparameterization $\boldsymbol{\zeta} = \mathbf{g}(\sigma^2, \boldsymbol{\xi}) \equiv (g_1(\sigma^2), \mathbf{g}_2(\boldsymbol{\xi}))$ chosen so that the new parameter space is unconstrained, that is, $\mathbb{R}^{r+1}$. For this alternative parameterization, as long as $\mathbf{g}$ is sufficiently regular, it is possible to calculate both the marginal ML estimate $\hat{\boldsymbol{\zeta}}$ and the associated expected information evaluated at the ML estimate $I_{\boldsymbol{\zeta}}(\hat{\boldsymbol{\zeta}})$, given that we can compute these quantities for the original parameterization. Consider the partition of $I_{\boldsymbol{\zeta}}(\hat{\boldsymbol{\zeta}})$ given by ($u$ is a scalar, $\mathbf{v}$ is an $r \times 1$ vector)

$$(4.2) \qquad I_{\boldsymbol{\zeta}}(\hat{\boldsymbol{\zeta}}) = \begin{bmatrix} u & \mathbf{v}' \\ \mathbf{v} & \mathbf{A} \end{bmatrix}$$

and define $\hat{\mathbf{V}} = \mathbf{A} - \mathbf{v}\mathbf{v}'/u$. We perform a Metropolis step in the alternative parameterization with a $t$-density proposal centered at the previous iteration, with $d$ degrees of freedom and scale matrix $c^2\hat{\mathbf{V}}$. To implement this MCMC method one only has to specify $d$ and $c$. The number of degrees of freedom is not an issue; one can even use a Gaussian instead of a $t$. As for $c$, one can always experiment with a few values, and a reasonable starting guess is $c \approx 2.4/\sqrt{r}$ [Gelman, Carlin, Stern and Rubin (1995)].

We next give the precise details of the algorithm. Assuming that the chain is currently at state $(\boldsymbol{\theta}_{(\text{old})}, \sigma^2_{(\text{old})}, \boldsymbol{\xi}_{(\text{old})})$, the next state $(\boldsymbol{\theta}_{(\text{new})}, \sigma^2_{(\text{new})}, \boldsymbol{\xi}_{(\text{new})})$ is determined by the following steps:

*Step* 1. Draw $\boldsymbol{\theta}_{(\text{new})} \sim \mathsf{N}((\mathbf{X}'\boldsymbol{\Sigma}^{-1}\mathbf{X})^{-1}\mathbf{X}'\boldsymbol{\Sigma}^{-1}\mathbf{y}, \sigma^2_{(\text{old})}(\mathbf{X}'\boldsymbol{\Sigma}^{-1}\mathbf{X})^{-1})$, where $\boldsymbol{\Sigma} \equiv \boldsymbol{\Sigma}(\boldsymbol{\xi}_{(\text{old})})$.

*Step* 2. Draw $\sigma^2_{(\text{new})} \sim \Gamma^{-1}(n/2 + a - 1, \mathbf{y}'\boldsymbol{\Sigma}^{-1}\mathbf{y}/2 + r_{\sigma^2})$, where $\boldsymbol{\Sigma} \equiv \boldsymbol{\Sigma}(\boldsymbol{\xi}_{(\text{old})})$ and by convention $a = 2$ in the case of the empirical Bayes method. Recall that in the case of the reference and independence Jeffreys priors $a = 1$, whereas in the Jeffreys-rule $a = 1 + q/2$. Also, $r_{\sigma^2}$ is the rate of the exponential prior in the empirical Bayes case and should be set to zero otherwise.

*Step* 3. Draw $\boldsymbol{\delta} \sim t_{(d)}(g_2(\boldsymbol{\xi}_{(\text{old})}), c^2\hat{\mathbf{V}})$ and let $\boldsymbol{\xi}_\star = \mathbf{g}_2^{-1}(\boldsymbol{\delta})$; then

$$\boldsymbol{\xi}_{(\text{new})} = \begin{cases} \boldsymbol{\xi}_\star, & \text{w.p. } \rho, \\ \boldsymbol{\xi}_{(\text{old})}, & \text{w.p. } 1 - \rho, \end{cases}$$

where

$$\rho = \min\left\{1, \frac{\mathsf{N}(\mathbf{y}|\mathbf{X}\boldsymbol{\theta}_{(\text{new})}, \sigma^2_{(\text{new})}\boldsymbol{\Sigma}(\boldsymbol{\xi}_{(\text{new})}))\pi(\boldsymbol{\xi}_{(\text{new})})}{\mathsf{N}(\mathbf{y}|\mathbf{X}\boldsymbol{\theta}_{(\text{new})}, \sigma^2_{(\text{new})}\boldsymbol{\Sigma}(\boldsymbol{\xi}_{(\text{old})}))\pi(\boldsymbol{\xi}_{(\text{old})})} \frac{q(\boldsymbol{\xi}_{(\text{old})}|\boldsymbol{\xi}_{(\text{new})})}{q(\boldsymbol{\xi}_{(\text{new})}|\boldsymbol{\xi}_{(\text{old})})}\right\}$$

and

$$q(x|y) = t_{(d)}(\mathbf{g}_2(x)|\mathbf{g}_2(y), c^2\hat{\mathbf{V}})\frac{\partial g_2}{\partial \boldsymbol{\xi}}\bigg|_{\boldsymbol{\xi}=x}.$$

4.5. *An example.* For illustrative purposes, we are going to consider a simplified version of a problem analyzed in detail in Bayarri et al. (2002). In that paper, the authors consider the analysis of a computer model that simulates the crash of prototype vehicles against a barrier, recording the velocity curve from the point of impact until the vehicle stops. For our purposes, we will imagine that if we input an impact velocity and an instant in time, the computer model will return the velocity of the vehicle at that point in time after impact.

The data consist of the output of the computer model at 19 values of $t$ and 9 initial velocities $v$, which corresponds to a 171-point design set that follows a Cartesian product. If we subtract from each of the individual curves the initial velocity, what happens is that all the curves will start at zero and decay at a rate roughly proportional to the initial velocity. For this transformed data, which is plotted in Figure 1, it seems reasonable to consider the mean function $\mathbb{E}Y(v,t) = vt\theta$, which in the notation of Section 2 translates into $q=1$ and $\psi(v,t) = vt$. In particular, this problem satisfies Assumptions $\mathcal{A}2$ and $\mathcal{A}3$.

In terms of the correlation function, we will assume a two-dimensional separable power exponential function, with the roughness parameters fixed at 2, and the parameterization

(4.3) $\quad c((t_1, v_1), (t_2, v_2)) = \exp(-\beta_t |t_1 - t_2|^2) \exp(-\beta_v |v_1 - v_2|^2).$

Figure 2 shows estimates of the posterior distribution associated with each of the priors that we have introduced in the paper. The empirical Bayes corresponds to centering the exponential priors at 10 times the computed ML estimates, and the dotted lines on the figure correspond to those densities. The sampling mechanism was implemented exactly as detailed above, and we used a $t$ density with three degrees of freedom and $c = 1.7$. The transformation **g** that we used was the logarithmic one. The acceptance rate was roughly 0.27 in all cases, and the results correspond to 50,000 iterations, minus 500 that were discarded as burn-in. The chains seem to reach stationarity very quickly, and 50,000 iterations is certainly more than we actually need for reliable inference.

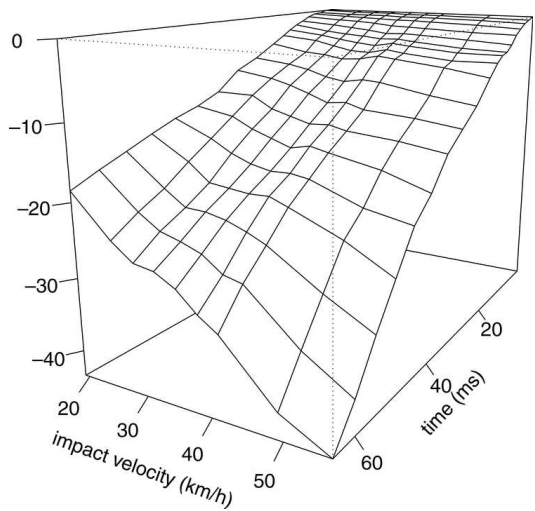

FIG. 1. *Relative velocity as a function of time after impact and of impact velocity.*

For these data the answers are quite robust with respect to the type of default prior chosen. In the next section we study finer details of the methods proposed in the paper, in particular by studying the frequentist properties of the resulting Bayesian procedures.

**5. Comparison of the priors.** When faced with more than one valid default prior specification strategy, it is often argued that one way to distinguish among these is by studying the frequentist properties of the resulting Bayesian inferential procedures. This usually takes the form of computing the frequentist coverage of the $(1-\gamma) \times 100\%$ equal-tailed Bayesian credible interval for one of the parameters of interest. The closer to the nominal level is this frequentist coverage, the "better" is the prior. For discussion and further references, see Kass and Wasserman (1996).

Our study was conducted in the context of the power exponential correlation function with $p=2$, parameterized as in (4.3), with the roughness parameters ($\alpha_1$ and $\alpha_2$) treated as known. For a $5 \times 5$ equally spaced grid in $[0,1] \times [0,1]$, we considered two options for the mean structure: an unknown constant term that we fixed at $\theta = 1$ and $\mathbb{E}[Y(\mathbf{x})|\boldsymbol{\eta}] = \theta x_1$, where again $\theta$ was fixed at 1. Several choices for the other parameters were considered, and we simulated 3000 draws from the ensuing model. Each Markov chain consisted of 15,000 iterations, the first 100 being discarded as burn-in. We calculated 95% credible intervals for $\sigma^2$, $\theta$ and $\beta_1$, and partial results are summarized in Table 2. Along with the estimate of the coverage probability, we present also an estimate of the expected length of the resulting credible interval and the standard deviation associated with the estimate.

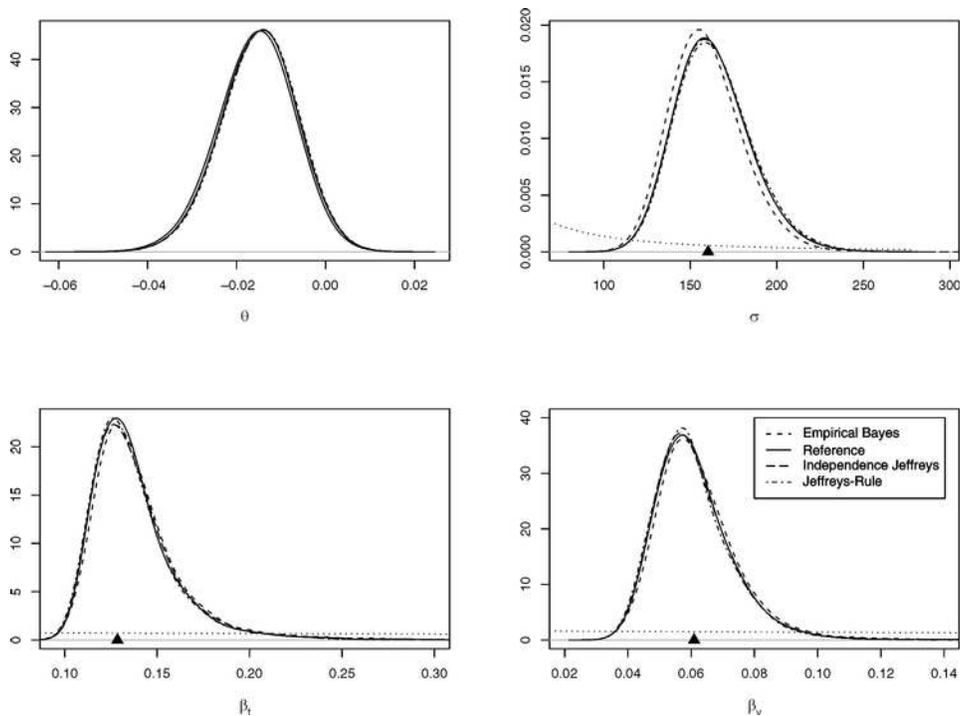

FIG. 2. *Smoothed histograms of samples drawn from the posterior distribution associated with each of the priors. Dotted lines correspond to exponential priors of empirical Bayes, and the triangles indicate the marginal likelihood estimates.*

A clear conclusion to extract from the results of these experiments is the comparatively poor behavior of the empirical Bayes method, which resulted from centering the priors at 10 times the computed ML estimate. This conclusion is particularly important since, as we have already pointed out, whenever a formal objective prior is not available, practitioners tend to resort to similar strategies to produce so-called "vague" or "diffuse" priors. On the basis of this study, one might argue against this type of approach to objective Bayesian analysis.

It is also possible to argue against the use of the Jeffreys-rule prior, and the reason for its inferior behavior has to do with the power $a = 1 + q/2$ in its formulation. This was already reported in Berger, De Oliveira and Sansó (2001): this type of prior is known to add spurious degrees of freedom to uncertainty statements.

The present study is not quite decisive about how the reference prior and the independence Jeffreys prior compare, and one may conclude that they present virtually equivalent performance. On the basis that the reference prior is computationally more demanding, we recommend the independence Jeffreys prior as a default prior for the problem at hand.

## APPENDIX

**A.0.1. Auxiliary facts.** Throughout the Appendix we will repeatedly use the following results. They are standard propositions whose proof can easily be found [e.g., Harville (1997) and Tong (1990)].

FACT 1. If for $i = 1, \ldots, r$ the product $\mathbf{A}_i \mathbf{B}_i$ is possible, then we have
$$\left(\bigotimes_{i=1}^{r} \mathbf{A}_i\right)\left(\bigotimes_{i=1}^{r} \mathbf{B}_i\right) = \bigotimes_{i=1}^{r} (\mathbf{A}_i \mathbf{B}_i).$$

FACT 2. If the matrices $\mathbf{A}_i$, $i = 1, \ldots, p$, are invertible, then $\bigotimes_{i=1}^{r} \mathbf{A}_i$ is invertible and one has
$$\left(\bigotimes_{i=1}^{r} \mathbf{A}_i\right)^{-1} = \bigotimes_{i=1}^{r} \mathbf{A}_i^{-1}.$$

TABLE 2
*Coverage probability of 95% credible intervals*

|  | **Coverage prob.** | **Expected length** | **Std. dev.** |
|---|---|---|---|
| Interval for $\sigma^2$; parameter vector fixed at $(1.5, 3.2, 3.6, 1.5, 1.7)$; $\mathbb{E}Y = 1$. | | | |
| Empirical Bayes | 0.857 | 2.904 | 0.005 |
| Reference prior | 0.956 | 11.354 | 0.003 |
| Independence Jeffreys prior | 0.954 | 11.698 | 0.003 |
| Jeffreys-rule prior | 0.946 | 6.334 | 0.003 |
| Interval for $\theta$; parameter vector fixed at $(1.5, 3.2, 3.6, 1.5, 1.7)$; $\mathbb{E}Y = \theta \equiv 1$. | | | |
| Empirical Bayes | 0.805 | 3.156 | 0.007 |
| Reference prior | 0.952 | 7.434 | 0.005 |
| Independence Jeffreys prior | 0.956 | 7.410 | 0.005 |
| Jeffreys-rule prior | 0.902 | 5.265 | 0.012 |
| Interval for $\theta$; parameter vector fixed at $(1.5, 0.2, 0.6, 1.5, 1.7)$; $\mathbb{E}Y = \theta x_1 \equiv x_1$. | | | |
| Empirical Bayes | 0.840 | 1.077 | 0.007 |
| Reference prior | 0.948 | 0.950 | 0.003 |
| Independence Jeffreys prior | 0.955 | 0.931 | 0.004 |
| Jeffreys-rule prior | 0.928 | 1.027 | 0.005 |
| Interval for $\beta_1$; parameter vector fixed at $(1.5, 0.2, 0.6, 1.5, 1.7)$; $\mathbb{E}Y = 1$. | | | |
| Empirical Bayes | 0.826 | 1.145 | 0.007 |
| Reference prior | 0.946 | 0.901 | 0.004 |
| Independence Jeffreys prior | 0.954 | 0.890 | 0.004 |
| Jeffreys-rule prior | 0.919 | 1.059 | 0.005 |
| Interval for $\beta_1$; parameter vector fixed at $(1.5, 0.2, 0.6, 1.0, 1.2)$; $\mathbb{E}Y = 1$. | | | |
| Empirical Bayes | 0.797 | 1.815 | 0.007 |
| Reference prior | 0.948 | 1.042 | 0.004 |
| Independence Jeffreys prior | 0.948 | 0.996 | 0.004 |
| Jeffreys-rule prior | 0.916 | 1.255 | 0.005 |

The order of the fixed parameter vector is $(\sigma^2, \beta_1, \beta_2, \alpha_1, \alpha_2)$.

FACT 3. *If the matrices* $\mathbf{A}_i$, $i = 1, \ldots, p$, *are of dimension* $n_i \times n_i$, *then with* $n_{(i)} = \prod_{k \neq i} n_k$,

$$\left| \bigotimes_{i=1}^{r} \mathbf{A}_i \right| = \prod_{i=1}^{r} |\mathbf{A}_i|^{n_{(i)}}.$$

FACT 4. *If* $\mathbf{A}$ *and* $\mathbf{B}$ *are* $m \times n$ *matrices, and* $\mathbf{C}$ *is* $p \times q$, *one has*

$$\mathbf{C} \otimes (\mathbf{A} + \mathbf{B}) = \mathbf{C} \otimes \mathbf{A} + \mathbf{C} \otimes \mathbf{B}, \qquad (\mathbf{A} + \mathbf{B}) \otimes \mathbf{C} = \mathbf{A} \otimes \mathbf{C} + \mathbf{B} \otimes \mathbf{C}.$$

FACT 5. $\operatorname{tr}(\mathbf{A} \otimes \mathbf{B}) = (\operatorname{tr} \mathbf{A})(\operatorname{tr} \mathbf{B})$.

FACT 6. $(\mathbf{A} \otimes \mathbf{B})' = \mathbf{A}' \otimes \mathbf{B}'$.

FACT 7. *If* $\mathbf{A}$ *is a function of* $\theta$ *but* $\mathbf{B}$ *is not, then*

$$\frac{\partial}{\partial \theta}[\mathbf{A} \otimes \mathbf{B}] = \left[\frac{\partial}{\partial \theta}\mathbf{A}\right] \otimes \mathbf{B}, \qquad \frac{\partial}{\partial \theta}[\mathbf{B} \otimes \mathbf{A}] = \mathbf{B} \otimes \left[\frac{\partial}{\partial \theta}\mathbf{A}\right].$$

FACT 8. *If* $\boldsymbol{\Sigma} \equiv \boldsymbol{\Sigma}(\theta)$ *is a positive definite matrix whose entries are differentiable with respect to* $\theta$, *then with* $\dot{\boldsymbol{\Sigma}} = \frac{\partial}{\partial \theta}\boldsymbol{\Sigma}$,

$$\frac{\partial}{\partial \theta}\log|\boldsymbol{\Sigma}| = \mathrm{tr}[\boldsymbol{\Sigma}^{-1}\dot{\boldsymbol{\Sigma}}], \qquad \frac{\partial}{\partial \theta}\boldsymbol{\Sigma}^{-1} = -\boldsymbol{\Sigma}^{-1}\dot{\boldsymbol{\Sigma}}\boldsymbol{\Sigma}^{-1}.$$

FACT 9. *Let* $\mathbf{X} \sim \mathsf{N}(\boldsymbol{\mu}, \boldsymbol{\Sigma})$ *and let* $\mathbf{A}$ *and* $\mathbf{B}$ *be symmetric matrices. Then*

$$\mathbb{E}\mathbf{X}'\mathbf{A}\mathbf{X} = \mathrm{tr}\,\mathbf{A}\boldsymbol{\Sigma} + \boldsymbol{\mu}'\mathbf{A}\boldsymbol{\mu},$$

$$\mathrm{Cov}(\mathbf{X}'\mathbf{A}\mathbf{X}, \mathbf{X}'\mathbf{B}\mathbf{X}) = 2\,\mathrm{tr}\,\mathbf{A}\boldsymbol{\Sigma}\mathbf{B}\boldsymbol{\Sigma} + 4\boldsymbol{\mu}'\mathbf{A}\boldsymbol{\Sigma}\mathbf{B}\boldsymbol{\mu}.$$

FACT 10. *Let* $\mathbf{A}$ *be a nonsingular matrix. Then* $|\mathbf{A} + \mathbf{1}\mathbf{1}'| = |\mathbf{A}|(1 + \mathbf{1}'\mathbf{A}^{-1}\mathbf{1})$. *If, furthermore,* $\mathbf{1}'\mathbf{A}^{-1}\mathbf{1} \neq -1$, *then* $\mathbf{A} + \mathbf{1}\mathbf{1}'$ *is nonsingular and*

$$(\mathbf{A} + \mathbf{1}\mathbf{1}')^{-1} = \mathbf{A}^{-1} - \frac{\mathbf{A}^{-1}\mathbf{1}\mathbf{1}'\mathbf{A}^{-1}}{1 + \mathbf{1}'\mathbf{A}^{-1}\mathbf{1}}.$$

**A.0.2. Proofs of Section 2.**

PROOF OF PROPOSITION 2.1. Define the shorthand $\ell^I = \log L^I(\sigma^2, \boldsymbol{\xi}|\mathbf{y})$. Berger, De Oliveira and Sansó (2001) prove that $\partial \ell^I/\partial \sigma^2 = (S^2_{\boldsymbol{\xi}} - \mathbb{E}S^2_{\boldsymbol{\xi}})/(2\sigma^4)$ and $\partial \ell^I/\partial \xi_i = (R^i_{\boldsymbol{\xi}} - \mathbb{E}R^i_{\boldsymbol{\xi}})/(2\sigma^2)$, where $S^2_{\boldsymbol{\xi}}/\sigma^2 \sim \chi^2_{n-q}$ and $R^i_{\boldsymbol{\xi}}$ is a quadratic form on $\mathbf{P}\mathbf{Y} \sim \mathsf{N}(\mathbf{0}, \sigma^2\mathbf{P}\boldsymbol{\Sigma})$, associated with the matrix $\boldsymbol{\Sigma}^{-1}\dot{\boldsymbol{\Sigma}}^i\boldsymbol{\Sigma}^{-1}$. The result follows from Fact 9 and elementary properties of the determinant function. □

PROOF OF PROPOSITION 2.2. Suppose $\mathbf{Y}|\boldsymbol{\theta}, \boldsymbol{\xi} \sim \mathsf{N}(\mathbf{X}\boldsymbol{\theta}, \boldsymbol{\Sigma})$, where $\boldsymbol{\Sigma} \equiv \boldsymbol{\Sigma}(\boldsymbol{\xi})$, and let $\ell \equiv \log L(\boldsymbol{\theta}, \boldsymbol{\xi}|\mathbf{y})$ be the associated log-likelihood function. Standard results of matrix differentiation plus Facts 8 and 9 of Appendix A.0.1 allow one to conclude that

$$\partial \ell/\partial \boldsymbol{\theta} = (\mathbf{X}'\boldsymbol{\Sigma}^{-1}\mathbf{y} - \mathbb{E}\mathbf{X}'\boldsymbol{\Sigma}^{-1}\mathbf{y})/2$$

and $\partial \ell/\partial \xi_i = (R^i_{\boldsymbol{\xi}} - \mathbb{E}R^i_{\boldsymbol{\xi}})/2$, where $R^i_{\boldsymbol{\xi}}$ is a quadratic form on $(\mathbf{Y} - \mathbf{X}\boldsymbol{\theta}) \sim \mathsf{N}(\mathbf{0}, \boldsymbol{\Sigma})$ associated with the matrix $\boldsymbol{\Sigma}^{-1}\dot{\boldsymbol{\Sigma}}^i\boldsymbol{\Sigma}^{-1}$. The result follows easily from elementary properties of the determinant function. □

**A.0.3. Proofs of Section 3.**

A.0.4. *Proofs of Section* 3.2.

PROOF OF LEMMA 3.1. We start out by proving (3.3). Since $\boldsymbol{\Sigma} = \bigotimes_{k=1}^{r} \boldsymbol{\Sigma}_k$, from Fact 2 of Appendix A.0.1 it follows that $\boldsymbol{\Sigma}^{-1} = \bigotimes_{k=1}^{r} \boldsymbol{\Sigma}_k^{-1}$. Then by Fact 6 of Appendix A.0.1, repeated application of Fact 1 and by Assumption $\mathcal{A}2$, one has $\mathbf{X}'\boldsymbol{\Sigma}^{-1}\mathbf{X} = \bigotimes_{k=1}^{r} \mathbf{X}_k'\boldsymbol{\Sigma}_k^{-1}\mathbf{X}_k$, which in conjunction with Fact 3 establishes the result. This last expression, Fact 2 and repeated application of Fact 1 show (3.2). □

PROOF OF LEMMA 3.2. Formulas (3.9)–(3.11) are shown in Berger, De Oliveira and Sansó (2001) to follow essentially from (3.4) and from the fact that $\lim_{\xi \to 0} \nu(\xi) = 0$, a simple consequence of (3.4) and $\rho(0|\xi) = 0$. The assumption $\mathbf{1}_{n_k}' \mathbf{D}_k^{-1} \mathbf{1}_{n_k} \neq 0$ and (3.5) assure that the matrices are well defined and that the expansions are meaningful.

We next prove (3.12). For simplicity, we drop the subscript $k$ in the sequel. Suppose then that $\mathbf{X} = \mathbf{1}$, and write $\boldsymbol{\Sigma} = \mathbf{A} + \mathbf{11}'$. Using Fact 10 and simple manipulations, Berger, De Oliveira and Sansó (2001) prove, in a more general context, that $\mathbf{1}(\mathbf{1}'\boldsymbol{\Sigma}^{-1}\mathbf{1})^{-1}\mathbf{1}' = \mathbf{1}(\mathbf{1}'\mathbf{A}^{-1}\mathbf{1})^{-1}\mathbf{1}' + \mathbf{11}'$. Pre- and postmultiplying this last equation by $\boldsymbol{\Sigma}^{-1} = \mathbf{A}^{-1} - (\mathbf{A}^{-1}\mathbf{11}'\mathbf{A}^{-1})/(1+\mathbf{1}'\mathbf{A}^{-1}\mathbf{1})$ and simplifying yields, after some algebra,

$$\boldsymbol{\Sigma}^{-1}\mathbf{1}(\mathbf{1}'\boldsymbol{\Sigma}^{-1}\mathbf{1})^{-1}\mathbf{1}'\boldsymbol{\Sigma}^{-1} = \frac{\mathbf{A}^{-1}\mathbf{11}'\mathbf{A}^{-1}}{(1+\mathbf{1}'\mathbf{A}^{-1}\mathbf{1})\mathbf{1}'\mathbf{A}^{-1}\mathbf{1}}.$$

Substituting $\mathbf{A} = \nu(\xi)(\mathbf{D} + o(1))$ in the right-hand side yields the first part of (3.12). The second part follows directly from formula $\boldsymbol{\Phi}_k = \boldsymbol{\Sigma}_k^{-1}\mathbf{X}_k(\mathbf{X}_k' \times \boldsymbol{\Sigma}_k^{-1}\mathbf{X}_k)^{-1}\mathbf{X}_k'\boldsymbol{\Sigma}_k^{-1}$ and expansion (3.9). Equation (3.5) assures that $\mathbf{H}_k$ is well defined, and it is also easy to show that $\mathbf{H}_k \neq \mathbf{0}$. □

PROOF OF PROPOSITION 3.3. Formula (3.13) is a simple consequence of Lemma 3.1 and Fact 3 of Appendix A.0.1. The continuity of $L^I(\boldsymbol{\xi} \mid \mathbf{y})$ as a function of $\boldsymbol{\xi}$ follows from the continuity of $\rho$ as a function of the parameter and the product form of the correlation function. We will use throughout this proof the facts that, as $\xi \to 0$, $\nu(\xi) \to 0$, and that, as $\xi_k \to \infty$, $\boldsymbol{\Sigma}_k \to \mathbf{I}_{n_k}$. These are simple consequences of the assumptions of Lemma 3.2.

Using (3.13) and the expansions of Lemma 3.2, it is possible to check that in the circumstances of part (a) of Proposition 3.3 we have

$$L^I(\boldsymbol{\xi}|\mathbf{y}) \propto \prod_{k \in B} \nu(\xi_k)^{(n_k-3+2a)/2} \prod_{k \in A \cap B} \nu(\xi_k)^{1/2}(1+g(\boldsymbol{\xi}))$$
$$\times \left\{ \mathbf{y}' \left[ \mathbf{C} - \prod_{k \in \bar{A} \cap B} \nu(\xi_k) \mathbf{C}^\star \right] \mathbf{y} \right\}^{-(n-3+2a)/2},$$

where $\mathbf{C} \equiv \bigotimes_{k=1}^{r} \mathbf{C}_k$, $\mathbf{C}^\star \equiv \bigotimes_{k=1}^{r} \mathbf{C}_k^\star$, with

$$\mathbf{C}_k = \begin{cases} \mathbf{G}_k, & \text{if } k \in B, \\ \mathbf{I}_{n_k}, & \text{if } k \in \bar{B}, \end{cases}$$

and

$$\mathbf{C}_k^\star = \begin{cases} \mathbf{F}_k, & \text{if } k \in \bar{A} \cap B, \\ \mathbf{H}_k, & \text{if } k \in A \cap B, \\ \mathbf{X}_k(\mathbf{X}_k'\mathbf{X}_k)^{-1}\mathbf{X}_k', & \text{if } k \in \bar{B}. \end{cases}$$

We next argue that the quantity within the braces can be bounded between two positive constants as long as $\xi_k$, $k \in B$, are sufficiently small, in which case the result follows immediately.

If $\bar{A} \neq \varnothing$, then that quantity converges to $\mathbf{y}'\mathbf{C}\mathbf{y}$, that we claim to be positive. To see that, recall that the Kronecker product of positive definite matrices is still a positive definite matrix. As a consequence, it suffices to show that $\mathbf{G}_k$ is positive definite. This follows from (3.9) and the fact that $\mathbf{\Sigma}_k^{-1}$ is positive definite (recall that $\nu$ is positive).

If $\bar{A} = \varnothing$, then it suffices to show that $\mathbf{y}'(\mathbf{C} - \mathbf{C}^\star)\mathbf{y} > 0$. To see that, note that it can be checked that $\mathbf{C} - \mathbf{C}^\star = \mathbf{C}[\mathbf{I}_n - \mathbf{X}(\mathbf{X}'\mathbf{C}\mathbf{X})^{-1}\mathbf{X}'\mathbf{C}]$. The quantity between brackets in the last expression is a projection matrix, and hence [cf. Harville (1997), page 262]

$$\mathbf{C}[\mathbf{I}_n - \mathbf{X}(\mathbf{X}'\mathbf{C}\mathbf{X})^{-1}\mathbf{X}'\mathbf{C}] = [\mathbf{I}_n - \mathbf{X}(\mathbf{X}'\mathbf{C}\mathbf{X})^{-1}\mathbf{X}'\mathbf{C}]'\mathbf{C}[\mathbf{I}_n - \mathbf{X}(\mathbf{X}'\mathbf{C}\mathbf{X})^{-1}\mathbf{X}'\mathbf{C}].$$

The result follows since, from what we have seen above, $\mathbf{C}$ is positive definite.

The proof of part (b) follows essentially from the same type of arguments. Let $\mathbf{C}$ and $\mathbf{C}^\star$ be as before but with $\mathbf{I}_{n_k}$ and $\mathbf{X}_k(\mathbf{X}_k'\mathbf{X}_k)^{-1}\mathbf{X}_k'$ replaced by $\mathbf{\Sigma}_k^{-1}$ and $\mathbf{\Phi}_k$, respectively. It is possible to check that

$$L^I(\boldsymbol{\xi}|\mathbf{y}) \propto \prod_{k \in B} \nu(\xi_k)^{(n_k - 3 + 2a)/2} \prod_{k \in A \cap B} \nu(\xi_k)^{1/2}(1 + g(\{\xi_k, k \in B\}))$$
$$\times \prod_{k \in \bar{B}} \{|\mathbf{\Sigma}_k|^{-n_{(k)}/2}(\mathbf{X}_k'\mathbf{\Sigma}_k^{-1}\mathbf{X}_k)^{-1/2}\}$$
$$\times \left\{\mathbf{y}'\left[\mathbf{C} - \prod_{k \in \bar{A} \cap B} \nu(\xi_k)\mathbf{C}^\star\right]\mathbf{y}\right\}^{-(n - 3 + 2a)/2}.$$

If $\bar{A} \neq \varnothing$, then the quantity within the braces in the last factor converges to $\mathbf{y}'\mathbf{C}\mathbf{y} > 0$. In this case, $h$ can be taken to be the second-to-last factor.

If $\bar{A} = \varnothing$, then we must show that $\mathbf{C} - \mathbf{C}^\star$ is positive definite in the set $\{\delta_k < \xi_k < M_k, k \in \bar{B}\}$. This follows from an argument similar to the one used for part (a), and therefore will be omitted. Here $h$ can be chosen to be the product of the last and second-to-last factors. $\square$

A.0.5. *Proofs of Section* 3.3.

PROOF OF PROPOSITION 3.4. Using Facts 1 and 7 of Appendix A.0.1 along with (3.2), it is clear that

$$\mathbf{W}_k = \mathbf{I}_{n_1} \otimes \cdots \otimes \dot{\mathbf{\Sigma}}_k \mathbf{\Sigma}_k^{-1} \otimes \cdots \otimes \mathbf{I}_{n_r} - \mathbf{\Sigma}_1 \mathbf{\Phi}_1 \otimes \cdots \otimes \dot{\mathbf{\Sigma}}_k \mathbf{\Phi}_k \otimes \cdots \otimes \mathbf{\Sigma}_r \mathbf{\Phi}_r.$$

Next, note that $\mathbf{\Sigma}_k \mathbf{\Phi}_k$ is a projection matrix, so that it is idempotent, and its rank, and therefore its trace, is $\operatorname{rank} \mathbf{X}_{n_k} = 1$. Also, $\mathbf{\Phi}_k \mathbf{\Sigma}_k \mathbf{\Phi}_k = \mathbf{\Phi}_k$.

These facts, Fact 5 of Appendix A.0.1, some algebra and the known fact that $\operatorname{tr} \mathbf{AB} = \operatorname{tr} \mathbf{BA}$ whenever the products are possible, allows one to show that

$$\operatorname{tr} \mathbf{W}_k = n_{(k)} \operatorname{tr} \dot{\mathbf{\Sigma}}_k \mathbf{\Sigma}_k^{-1} - \operatorname{tr} \dot{\mathbf{\Sigma}}_k \mathbf{\Phi}_k,$$

$$\operatorname{tr} \mathbf{W}_k^2 = n_{(k)} \operatorname{tr}(\dot{\mathbf{\Sigma}}_k \mathbf{\Sigma}_k^{-1})^2 + \operatorname{tr}(\dot{\mathbf{\Sigma}}_k \mathbf{\Phi}_k)^2 - 2 \operatorname{tr} \dot{\mathbf{\Sigma}}_k \mathbf{\Sigma}_k^{-1} \dot{\mathbf{\Sigma}}_k \mathbf{\Phi}_k,$$

$$\operatorname{tr} \mathbf{W}_i \mathbf{W}_j = \prod_{l \neq i,j} n_l \times \operatorname{tr} \dot{\mathbf{\Sigma}}_i \mathbf{\Sigma}_i^{-1} \operatorname{tr} \dot{\mathbf{\Sigma}}_j \mathbf{\Sigma}_j^{-1} - \operatorname{tr} \dot{\mathbf{\Sigma}}_i \mathbf{\Phi}_i \operatorname{tr} \dot{\mathbf{\Sigma}}_j \mathbf{\Phi}_j,$$

the key point to note being that only $\operatorname{tr} \mathbf{W}_i \mathbf{W}_j$ depends on more than one parameter. Therefore, we can write (2.7) as

$$(A.1) \qquad I_R(\boldsymbol{\xi}) = \begin{pmatrix} \mathbf{T} & \boldsymbol{\omega} \\ \boldsymbol{\omega}' & g_r^R(\xi_r) \end{pmatrix},$$

where $\mathbf{T}$ is $(p-1) \times (p-1)$ and does not depend on $\xi_r$, $\boldsymbol{\omega}$ is $(p-1) \times 1$ and $g_r^R(\xi_r) = \operatorname{tr} \mathbf{W}_r^2$, which depends on $\xi_r$ only. As a consequence of the above block format, we have [cf. Harville (1997), page 188]

$$|I_R(\boldsymbol{\xi})| = |\mathbf{T}|(g_r^R(\xi_r) - \boldsymbol{\omega}' \mathbf{T}^{-1} \boldsymbol{\omega}) \leq |\mathbf{T}| g_r^R(\xi_r),$$

the last step following from the fact that $\mathbf{T}^{-1}$ is positive definite. If we repeat this procedure $p$ times, we end up proving (3.18) by defining $\pi_k(\xi_k) = [g_k^R(\xi_k)]^{1/2} = [\operatorname{tr} \mathbf{W}_k^2]^{1/2}$.

We now study $\pi_k(\xi_k)$ as $\xi_k \to 0$. It is easy to verify that $\mathbf{D}_k \mathbf{G}_k$ is idempotent and that $\operatorname{tr} \mathbf{D}_k \mathbf{G}_k = n_k - 1$. Also, it is possible to show that $\mathbf{D}_k \mathbf{G}_k \mathbf{D}_k \mathbf{H}_k = \mathbf{D}_k \mathbf{H}_k$ and that $\operatorname{tr} \mathbf{D}_k \mathbf{H}_k = 1$. Note also that, as a consequence of the assumptions we have, as $\xi_k \to 0$,

$$(A.2) \qquad \dot{\mathbf{\Sigma}}_k = \nu'(\xi_k) \mathbf{D}_k (1 + o(1)).$$

These properties and simple expansions along with (3.9) and (3.12) show (3.19). Assumption $\mathcal{A}8$ is capital to this result.

Since as $\xi_k \to \infty$, $\mathbf{\Sigma}_k \to \mathbf{I}_{n_k}$, the integrability of $\pi_k(\xi_k)$ at infinity follows from Assumption $\mathcal{A}7$. □

PROOF OF PROPOSITION 3.5. This proof essentially mimics the one of Proposition 3.4. Using the definition of $\mathbf{U}_k$ and Facts 1 and 7 of Appendix A.0.1, it is easy to check that

$$\mathbf{U}_k = \mathbf{I}_{n_1} \otimes \cdots \otimes \dot{\mathbf{\Sigma}}_k \mathbf{\Sigma}_k^{-1} \otimes \cdots \otimes \mathbf{I}_{n_r}.$$

Next, essentially the same type of arguments as before allows one to write

$$\operatorname{tr} \mathbf{U}_k = n_{(k)} \operatorname{tr} \dot{\mathbf{\Sigma}}_k \mathbf{\Sigma}_k^{-1},$$
$$\operatorname{tr} \mathbf{U}_k^2 = n_{(k)} \operatorname{tr}(\dot{\mathbf{\Sigma}}_k \mathbf{\Sigma}_k^{-1})^2,$$
$$\operatorname{tr} \mathbf{U}_i \mathbf{U}_j = \prod_{l \neq i,j} n_l \times \operatorname{tr} \dot{\mathbf{\Sigma}}_i \mathbf{\Sigma}_i^{-1} \operatorname{tr} \dot{\mathbf{\Sigma}}_j \mathbf{\Sigma}_j^{-1}.$$

Again, only $\operatorname{tr} \mathbf{U}_i \mathbf{U}_j$ depends on more than one parameter, and the construction carried out in the proof of Proposition 3.4 can be repeated to show (3.20), where now $\pi_k^{J1}(\xi_k) \propto (\operatorname{tr} \mathbf{U}_k^2)^{1/2}$ and $\pi_k^{J2}(\xi_k) \propto \pi_k^{J1}(\xi_k) |\mathbf{X}_k' \mathbf{\Sigma}_k^{-1} \mathbf{X}_k|^{1/2}$.

To obtain the behavior of these functions as $\xi_k \to 0$, it suffices to use expansions (3.9), (3.11) and (A.2). It is also necessary to recall that $\mathbf{D}_k \mathbf{G}_k$ is idempotent and that its trace is $n_k - 1$ (and hence Assumption A8).

The integrability at infinity of $\pi_k^{Ji}(\cdot)$ follows from the fact that, as $\xi_k \to \infty$, $\mathbf{\Sigma}_k \to \mathbf{I}_{n_k}$, and Assumption A7. $\square$

A.0.6. *Proofs of Section* 3.4.

PROOF OF THEOREM 3.6. We will determine conditions under which

$$0 < \int_0^\infty L^I(\boldsymbol{\xi}|\mathbf{y}) \pi(\boldsymbol{\xi}) \, d\xi_1 \cdots d\xi_r < \infty$$

holds. To simplify the notation, we will write $f(\boldsymbol{\xi}) = L^I(\boldsymbol{\xi}|\mathbf{y})$. The function $\pi(\boldsymbol{\xi})$ will be associated either with the reference prior or with one of the Jeffreys priors we considered, but we will carry out the calculations for a general exponent $a$ in the integrated likelihood. In all of these instances of $\pi$, there are functions $\pi_k(\xi_k)$, $k = 1, \ldots, r$, such that $\pi(\boldsymbol{\xi}) \leq \prod_{k=1}^r \pi_k(\xi_k)$, and those have essentially the same behavior across the different choices for $\pi(\boldsymbol{\xi})$, so that a common proof can be sought—compare Propositions 3.4 and 3.5. Note that

$$(\text{A.3}) \qquad \int_0^\infty f(\boldsymbol{\xi}) \pi(\boldsymbol{\xi}) \, d\xi_1 \cdots d\xi_r \leq \int_0^\infty f(\boldsymbol{\xi}) \pi_1(\xi_1) \, d\xi_1 \cdots \pi_r(\xi_r) \, d\xi_r,$$

and therefore we will have posterior propriety if the right-hand side is finite.

The right-hand side of (A.3) can be partitioned into a sum of integrals, each of which has different regions of integration. If the region is of the form $[\delta_1, +\infty[ \times \cdots \times [\delta_r, +\infty[$, then there are no issues to address since the $\pi_k$ are integrable at infinity and in this set the integrated likelihood is bounded.

Therefore, we need only consider regions of the following type: let $B$ be a nonempty but otherwise arbitrary subset of $\{1,\ldots,n\}$, and define said regions by:

(i) for $k \in B$, $0 < \xi_k < \delta_k$, where $\delta_k$ can be chosen arbitrarily small; for $k \in \bar{B}$, $\xi_k > M_k$, where $M_k$ can be assumed arbitrarily large;

(ii) for $k \in B$, $0 < \xi_k < \delta_k$, where $\delta_k$ can be chosen arbitrarily small; for $k \in \bar{B}$, $\delta_k < \xi_k < M_k$, where $M_k$ is finite but otherwise arbitrary.

Let us consider (i) first. Combining expansions (3.14), (3.19) and (3.21), it is easy to see that

$$f(\boldsymbol{\xi}) \prod_{k=1}^{r} \pi_k(\xi_k) \propto \prod_{k \in A \cap B} |\nu'(\xi_k)| \nu(\xi_k)^{(n_k - 2 - 2\alpha_k + 2a)/2}$$
$$\times \prod_{k \in \bar{A} \cap B} |\nu'(\xi_k)| \nu(\xi_k)^{(n_k - 5 + 2a)/2}$$
$$\times \prod_{k \in \bar{B}} \pi_k(\xi_k)(1 + g(\boldsymbol{\xi})),$$

where $\alpha_k$ is defined in Proposition 3.5 ($\alpha_k \equiv 1$ in the case of the reference prior).

When we look at situation (ii), formula (3.15) and the same expansions as before for the priors yield an expression formally identical to the one above multiplied by $h(\xi_k, k \in \bar{B})$.

Since the functions $\pi_k$ are integrable at infinity, and for small enough $\xi$, $\nu(\xi) < 1$, it is clear that posterior propriety will be achieved for all three priors whenever the function $|\nu'(\xi)| \nu(\xi)^{(n_k - 5 + 2a)/2}$ is integrable at the origin for every $k$. This reduces to (3.23) since $\nu = o(1)$ and $\min n_k \geq 2$. $\square$

**Acknowledgment.** The author would like to thank Professor James O. Berger for his guidance and for numerous insightful discussions.

## REFERENCES

Bayarri, M. J., Berger, J. O., Higdon, D., Kennedy, M. C., Kottas, A., Paulo, R., Sacks, J., Cafeo, J. A., Cavendish, J., Lin, C. H. and Tu, J. (2002). A framework for validation of computer models. Technical Report 128, National Institute of Statistical Sciences.

Berger, J. O. and Bernardo, J. M. (1992). On the development of reference priors (with discussion). In *Bayesian Statistics 4* (J. M. Bernardo, J. O. Berger, A. P. Dawid and A. F. M. Smith, eds.) 35–60. Oxford Univ. Press. MR1380269

Berger, J. O., De Oliveira, V. and Sansó, B. (2001). Objective Bayesian analysis of spatially correlated data. *J. Amer. Statist. Assoc.* **96** 1361–1374. MR1946582

Bernardo, J. M. (1979). Reference posterior distributions for Bayesian inference (with discussion). *J. Roy. Statist. Soc. Ser. B* **41** 113–147. MR547240

CRESSIE, N. A. C. (1993). *Statistics for Spatial Data*, revised ed. Wiley, New York. MR1239641

FUENTES, M. (2003). Testing for separability of spatio-temporal covariance functions. Technical Report 2545, Dept. Statistics, North Carolina State Univ.

GELMAN, A., CARLIN, J. B., STERN, H. S. and RUBIN, D. B. (1995). *Bayesian Data Analysis*. Chapman and Hall, London. MR1385925

HARVILLE, D. A. (1974). Bayesian inference for variance components using only error contrasts. *Biometrika* **61** 383–385. MR368279

HARVILLE, D. A. (1997). *Matrix Algebra from a Statistician's Perspective*. Springer, Berlin. MR1467237

KASS, R. and WASSERMAN, L. (1996). The selection of prior distributions by formal rules. *J. Amer. Statist. Assoc.* **91** 1343–1370. MR1478684

KENNEDY, M. C. and O'HAGAN, A. (2000). Predicting the output from a complex computer code when fast approximations are available. *Biometrika* **87** 1–13. MR1766824

KENNEDY, M. C. and O'HAGAN, A. (2001). Bayesian calibration of computer models. *J. R. Stat. Soc. Ser. B Stat. Methodol.* **63** 425–464. MR1858398

SACKS, J., WELCH, W. J., MITCHELL, T. J. and WYNN, H. P. (1989). Design and analysis of computer experiments (with discussion). *Statist. Sci.* **4** 409–435. MR1041765

SHORT, M. and CARLIN, B. (2003). Multivariate spatiotemporal CDFs with measurement error. Technical Report 2003-015, Division of Biostatistics, Univ. Minnesota.

TONG, Y. L. (1990). *The Multivariate Normal Distribution*. Springer, Berlin. MR1029032

NATIONAL INSTITUTE OF STATISTICAL SCIENCES
STATISTICAL AND APPLIED MATHEMATICAL SCIENCES INSTITUTE
19 T. W. ALEXANDER DRIVE
P.O. BOX 14006
RESEARCH TRIANGLE PARK, NORTH CAROLINA 27709-4006
USA
E-MAIL: rui@niss.org